\documentclass{amsart}
\usepackage{amssymb}
\usepackage{amsfonts}
\usepackage{latexsym}

\newtheorem{theorem}{Theorem}[section]
\newtheorem{lemma}[theorem]{Lemma}
\newtheorem{proposition}[theorem]{Proposition}
\newtheorem{corollary}[theorem]{Corollary} 
\theoremstyle{definition}  
\newtheorem{definition}[theorem]{Definition}
\newtheorem{example}[theorem]{Example}

\theoremstyle{remark}
\newtheorem{remark}[theorem]{Remark}

\renewcommand{\text}[1]{\mbox{{\rm #1}}}
 
\newcommand{\id}{\text{id}}

\newcommand{\rank}{\text{rank\,}} 
 
\newcommand{\Map}{\text{Map}}
\newcommand{\End}{\text{End}} 

\newcommand{\Hom}{\text{Hom}} 
\newcommand{\Ad}{\text{Ad}}
\newcommand{\Sh}{\text{Sh}}
\newcommand{\Ind}{\text{Ind}}

\newcommand{\St}{\text{St}}

\newcommand{\comp}{\text{comp}}

\newcommand{\ti}{\tilde}   
\newcommand{\eps}{\varepsilon}

\newcommand{\I}{\mathcal{I}}
\renewcommand{\O}{\mathcal{O}}

\newcommand{\G}{\widehat{G}}
\newcommand{\K}{\widehat{K}}
\newcommand{\wG}{\widehat{G}}
\newcommand{\wB}{\widehat{B}}
\newcommand{\wH}{\widehat{H}}

\newcommand{\g}{\mathfrak{g}}
\newcommand{\h}{\mathfrak{h}}
\newcommand{\n}{\mathfrak{n}}

\newcommand{\wg}{\widehat{\mathfrak{g}}}
\newcommand{\wh}{\widehat{\mathfrak{h}}}

\newcommand{\wb}{\widehat{\mathfrak{b}}}

\newcommand{\la}{{<}}  %{\langle\,} 
\newcommand{\ra}{{>}}  %{\,\rangle} 

\begin{document}

\title{Dynamical twists in group algebras}  

\author{Pavel Etingof}
\address{Columbia University,
Department of Mathematics,
2990 Broadway, New York, NY 10027 and
MIT, Department of Mathematics,
77 Massachusetts Avenue, Cambridge, MA 02139-4307}

\author{Dmitri Nikshych} 
\address{MIT, Department of Mathematics,
77 Massachusetts Avenue 2-130, Cambridge, MA 02139-4307}
  
\date{January 29, 2001}

\begin{abstract} 
We classify dynamical twists in group algebras of finite groups.
Namely, we set up a bijective correspondence between gauge equivalence
classes of dynamical twists (which are solutions of a certain
non-linear functional equation) and isomorphism classes of 
``dynamical data'' described in purely group theoretical terms. 
This generalizes the classification
of usual twists obtained by Movshev and Etingof-Gelaki.
\end{abstract} 
\maketitle  

%%%%%%%%%%%%%%%%%%%%%%%%%%%%%%%%%%%%%%%%%%%%%%%%%%%%%%%%%%%%
%%%%%%%%%% INTRODUCTION         %%%%%%%%%%%%%%%%%%%%%%%%%%%%
%%%%%%%%%%%%%%%%%%%%%%%%%%%%%%%%%%%%%%%%%%%%%%%%%%%%%%%%%%%%

\begin{section}
{Introduction}In the last few years, following the pioneering paper \cite{F},
the theory of quantum groups has developed a new 
branch -- the theory of dynamical $R$-matrices and dynamical quantum groups
(see \cite{ES} for a review). 
In this theory, there is a useful notion of a dynamical twist, 
which was perhaps first introduced in \cite{BBB}. 

Namely, let $H$ be a Hopf algebra
over a field $k$, and $A$ be a finite Abelian group 
of group-like elements of $H$. 
A dynamical twist
is an $A$-invariant function  $A^*\to H\otimes H$, which satisfies a certain
nonlinear  functional equation (the so-called dynamical non-abelian 
$2$-cocycle condition), see Definition~\ref {dynamical twist}.

If $A$ is trivial, a dynamical twist
is just an ordinary twist, in the sense of Drinfeld (see, e.g., \cite{M}). 
Thus, the notion of a dynamical twist generalizes that of a usual twist. 

The significance of the notion of a dynamical twist 
consists in the fact that given  a dynamical twist $J(\lambda)$, 
one can endow the algebra $H\otimes \End_k k[A]$ 
with a nontrivial structure of a weak Hopf algebra (see \cite{EN}),
which is quasitriangular if so is $H$.
This is done using ``twisting'' of $H$ by $J$, 
which is a procedure analogous to 
usual twisting of Hopf algebras.  
In particular, for a quasitriangular 
$H$ with the $R$-matrix $\mathcal R$, the function 
${\mathcal R}(\lambda)=J^{-1}(\lambda)^{21}{\mathcal R}J(\lambda)$ 
is a dynamical $R$-matrix, i.e., it satisfies the quantum dynamical 
Yang-Baxter  equation (see \cite{ES}). 

In this paper we classify dynamical twists in group algebras 
of finite groups over an algebraically closed field $k$ of 
characteristic zero. Since group algebras are trivially quasitriangular, 
all such twists give rise to dynamical $R$-matrices. 
Our classification generalizes the classification 
of usual twists in group algebras, which was done in \cite{EG}, 
following closely the paper \cite{Mo}. 

Namely, let $A$ be an abelian  subgroup of a finite group $G$.
A dynamical datum for $(G, A)$ is a subgroup $K$ of $G$ together with 
a family of irreducible projective representations of $K$ satisfying a certain 
coherence condition, see Definition~\ref{dynamical datum}.
Our main result is the following 
\smallskip

\textbf{Theorem~\ref{our lovely correspondence}}. 
There is a bijection between
\newline
(i) gauge equivalence classes of dynamical twists, 
$J : A^*\to k[G]\otimes k[G]$, 
\newline
(ii) isomorphism classes of dynamical data for $(G, A)$.
\smallskip

The structure of the paper is as follows.

In Section $2$  we recall the definitions of a dynamical twist
and dynamical gauge equivalence.

In Section $3$ we use the idea of Movshev \cite{Mo} to associate
a family of semisimple algebras and bimodules
with every dynamical twist. 
Here we also define the notions of minimal and minimizable dynamical twists.

The concept of a dynamical datum for a pair $(G, A)$ 
is introduced in Section~$4$.

The main result of Section $5$ is Theorem~\ref{twist to datum}
which shows that every dynamical twist gives rise to an
isomorphism class of dynamical data
in such a manner that gauge equivalent twists give the same class of data.

The converse to this theorem is provided
in Section $6$, where we employ the exchange construction of \cite{EV} to 
construct a gauge equivalence class of
dynamical twists from a dynamical datum. The correspondence
between classes of twists and classes of data
is shown to be bijective in Theorem~\ref{our lovely correspondence}
which is the central result of this paper.  
Next, we explicitly
construct an example of a dynamical datum leading to a non-minimizable 
dynamical twist (Example~\ref{non-minimizable data}) --
which is a purely ``dynamical'' phenomenon, that does not exist for 
usual twists. We also explicitly
compute a family of dynamical twists in the group algebra of the
group $F_p^\times  \ltimes F_p$ of affine transformations of
the line over the field $F_p$ of $p$ elements
(Example~\ref{affine group}).

Finally, in Section $7$ we present a construction of minimal
dynamical twists in group algebras of 
finite nilpotent groups. Our technique here uses
the exponential map for nilpotent Lie algebras over finite fields
and a result of Kazhdan \cite{K}.

We expect that the methods of this paper can be applied to 
constructing and classifying
dynamical twists in universal enveloping algebras and, in particular,
to reproving and possibly improving the main result of \cite{X},
stating that a splittable triangular dynamical $r$-matrix can be quantized.
This is a subject of our future research. 

\smallskip

\textit{Acknowledgments}. 
We are grateful to R.~Guralnick for providing us with 
Example~\ref{non-minimizable data}. The authors were partially
supported by the NSF grant DMS-9988796.
The first author's work was partially done for the Clay
Mathematics Institute.
The second author thanks MIT for the hospitality during his visit.

\end{section}

%%%%%%%%%%%%%%%%%%%%%%%%%%%%%%%%%%%%%%%%%%%%%%%%%%%%%%%%%%%%
%%%%%%%%%%  DYNAMICAL TWISTS    %%%%%%%%%%%%%%%%%%%%%%%%%%%%
%%%%%%%%%%%%%%%%%%%%%%%%%%%%%%%%%%%%%%%%%%%%%%%%%%%%%%%%%%%%

\begin{section}
{Dynamical twists in Hopf algebras}

In this section we recall the definition of a dynamical twist \cite{BBB}.

Let $H$ be a Hopf algebra and $A$
be a finite Abelian subgroup of the group of group-like elements of $H$.
Then $k[A] = \Map({A^*} , k)$ is a commutative and cocommutative
Hopf subalgebra of $H$.  Let $P_\mu,\, \mu\in  {A^*}$,
be the minimal idempotents in $k[A]$ corresponding to characters of $A$.

\begin{definition}
\label{weight subspace}
For any left $H$-module $X$ and character $\mu\in  {A^*}$ define a
{\em weight subspace}
\begin{equation}
X[\mu] =\{x\in X \mid ax=\mu(a)x,\, \text{ for all } a\in A \}.
\end{equation}
\end{definition}

\begin{definition}
\label{zero weight}
We say that an element $x\in H^{\otimes n},\, n\geq 1$ 
has {\em zero weight} if $x\in H^{\otimes n}[0]$, where
$ H^{\otimes n}$ is viewed as a left $H$-module via the adjoint action;
in other words, if $x$ commutes with 
$\Delta^n(a)$ for all $a\in A$, where $\Delta^n: k[A] \to k[A]^{\otimes n}$ 
is the iterated comultiplication.
\end{definition}

\begin{definition}
\label{dynamical twist}
Let $J(\lambda) : A^* \to H\otimes H$ be a
zero weight function with invertible values. We say that  $J(\lambda)$ is a 
{\em dynamical twist}  in $H$ if it satisfies the following 
functional equations 
\begin{eqnarray}
\label{dynamical equation}
J^{12,3}(\lambda) J^{12}(\lambda-h^{(3)})
&=& J^{1,23}(\lambda) J^{23}(\lambda), \\
\label{counit of J}
(\eps\otimes \id)J(\lambda) = (\id \otimes \eps)J(\lambda)  &=& 1.
\end{eqnarray}
Here $J^{12,3}(\lambda) =(\Delta\otimes\id)J(\lambda)$,
$J^{12}(\lambda) = J(\lambda)\otimes 1$ etc.
The notation $\lambda\pm h^{(i)}$ means that
the argument $\lambda$ is shifted by the weight of the $i$-th component,
e.g., $J(\lambda-h^{(3)})$ is the element of $H\otimes H \otimes k[A]$
such that $\mu_3 (J(\lambda-h^{(3)})) = J(\lambda-\mu)$, where 
the index $i$ in $\mu_i$ indicates that $\mu$ is applied 
to the $i$th component.
\end{definition}

\begin{remark}
\label{versions of J}
There are two other versions of the dynamical twist 
condition~(\ref{dynamical equation}) that appeared in the literature :
\begin{eqnarray}
\label{+}
J^{12,3}(\lambda) J^{12}(\lambda+h^{(3)})
&=& J^{1,23}(\lambda) J^{23}(\lambda), \\
\label{+-}
J^{12,3}(\lambda) J^{12}(\lambda+h^{(3)})
&=& J^{1,23}(\lambda) J^{23}(\lambda-h^{(1)}).
\end{eqnarray}
Note that (\ref{dynamical equation}) is obtained from equations
(\ref{+}) and (\ref{+-}) by the changes of variable 
$\lambda\mapsto -\lambda$ and  
$\lambda\mapsto \lambda+h^{(1)}+h^{(2)}$ respectively.
\end{remark}

\begin{definition}
\label{dynamical gauge}
If $J(\lambda)$ is a dynamical twist in $H$ and $t(\lambda) : A^*\to H$ is a
zero weight function with invertible values
such that $\eps(t(\lambda))\equiv 1$, then 
\begin{equation}
J^t(\lambda) = \Delta(t(\lambda)^{-1})\,J(\lambda)\, 
                (t(\lambda-h^{(2)})\otimes t(\lambda))
\end{equation}
is also a dynamical twist in $H$, {\em gauge equivalent} to $J(\lambda)$.
The function    $t(\lambda)$ is called a {\em gauge transformation}.
\end{definition}

\begin{remark}
\label{single value}
Note (\cite{ES}, Appendix C) that a dynamical twist $J(\lambda)$
on $H$ is completely defined by its value $J(\lambda_0)$,
at any point $\lambda_0\in A^*$. 
Indeed, from Equation~(\ref{dynamical equation}) we have
\begin{equation}
J(\lambda_0-\mu) =
\mu_3  (J^{12,3}(\lambda_0)^{-1} J^{1,23}(\lambda_0) J^{23}(\lambda_0)), 
\end{equation}
for all $\mu\in A^*$.
\end{remark}

\begin{remark}
\label{J defines J(lambda)}
Let $H=k[G]$ be a group Hopf algebra and  $J\in H\otimes H$ be 
an invertible element with the properties
$(\eps\otimes \id)J = (\id \otimes \eps)J  = 1$
and
\begin{equation*}
(J^{12,3})^{-1} J^{1,23} J^{23} \in k[G]\otimes k[G]\otimes k[A].
\end{equation*}
Then one can check by a direct computation (\cite{ES}, 12.1) that
\begin{equation}
J(\lambda) = -\lambda_3( (J^{12,3})^{-1} J^{1,23} J^{23}).
\end{equation}
is a dynamical twists in $H$.
\end{remark}

\begin{remark}
If $A=\{1\}$ then the definition of dynamical twist
coincides with the usual notion of twist introduced by Drinfeld.
\end{remark}

\begin{example}
If $H=k[A]$, then $H$ is commutative and can be identified with the algebra
of functions on $A^*$.  Let $P_\mu,\, \mu\in A^*$ be  the minimal idempotents 
of $H$. Let $c: A^*\times A^* \to k^\times$ be any function
such that $c(\lambda,\, 0) = c(0,\,\lambda) = 1$. 
Then  $J =\sum_{\mu,\nu}\, c(\mu,\, \nu) P_\mu \otimes P_\nu$
satisfies the conditions of Remark~\ref{J defines J(lambda)} and hence
\begin{equation}
J(\lambda) = \sum_{\mu,\nu}\,
c(\mu+\nu,\,-\lambda)^{-1}\,c(\mu,\,\nu-\lambda)\,c(\nu,\,-\lambda)
P_\mu \otimes P_\nu.
\end{equation}
is a dynamical twist. Furthermore, every dynamical twist in $k[A]$
is of this form by Remarks~\ref{single value} and \ref{J defines J(lambda)}. 

In fact, it turns out that in this case $J(\lambda)$ 
is always gauge equivalent to the constant twist $1\otimes 1$.
Namely, consider a gauge transformation 
\begin{equation}
t(\lambda) = \sum_\mu\, c(\mu, -\lambda) P_\mu,\qquad \lambda\in A^*, 
\end{equation}
then $J(\lambda)=\Delta(t(\lambda)^{-1})
(t(\lambda-h^{(2)})\otimes t(\lambda))$.
%i.e., it is gauge equivalent to the constant twist $J(\lambda)=1\otimes 1$.
\end{example}

\end{section}

%%%%%%%%%%%%%%%%%%%%%%%%%%%%%%%%%%%%%%%%%%%%%%%%%%%%%%%%%%%%%%%%%%%%%%%%%%
%%%%%%%%%%%%%%  Algebras associated with a twist                   %%%%%%%
%%%%%%%%%%%%%%%%%%%%%%%%%%%%%%%%%%%%%%%%%%%%%%%%%%%%%%%%%%%%%%%%%%%%%%%%%%

\begin{section}
{Algebras $B_\lambda$ associated with a dynamical twist}
\label{algebras}

From now on let $H =k[G]$, the Hopf algebra of a finite group $G$,
then $A$ is a subgroup of $G$. Define $k[G/A]$ to be the quotient
of $k[G]$ by the left ideal generated by the elements $(a-1), \, a\in A$.

For any $\lambda\in {A^*}$ let us define a comultiplication and
counit on $k[G/ A]$ as follows 
\begin{equation*}
\Delta_\lambda(g) = (g \otimes g)J(\lambda), \qquad \eps_\lambda(g) = 1,
\end{equation*}
for all cosets $g\in G/ A$. Note that since $J(\lambda)$
has zero weight, the above operations are well-defined.

\begin{proposition}
\label{coalgebra c}
$C_\lambda = (k[G/ A], \Delta_\lambda, \eps_\lambda)$ is
a coassociative $G$-coalgebra (where $G$ acts via the left
multiplication) with counit.
\end{proposition}
\begin{proof}
It follows from the dynamical twist identity (\ref{dynamical equation}) that
\begin{equation}
\label{twist mod A}
(\Delta\otimes \id)J(\lambda) (J(\lambda) \otimes 1)
= (\id \otimes\Delta) J(\lambda) (1\otimes J(\lambda) )
\,\, \text{in} \,\, k[G]\otimes k[G]\otimes k[G/ A]
\end{equation}
which implies coassociativity
of $\Delta_\lambda$. The counit axiom is obvious.
\end{proof}

Let $B_\lambda$ be the associative algebra dual to $C_\lambda$.
Then $B_\lambda$ can be naturally identified with the algebra $F_0[G]$ 
of all functions on $G$ invariant under right translations by elements 
of $A$ with multiplication given by
\begin{equation}
f_1 *_\lambda f_2 (g) = f_1 \otimes f_2 ((g\otimes g)J(\lambda)),
\end{equation}
for all $f_1, f_2 \in F_0(G)$ and $g\in G$. It is a left $G$-algebra
with the action of $G$ given by
\begin{equation}
\label{action on FG}
h\circ f (g) = f(h^{-1}g),
\end{equation}
for all $f\in F_0[G]$ and $h,g\in G$.

\begin{remark}
\label{ge --> iso of Bl}
If $t(\lambda)$ is a gauge transformation  then coalgebras $C_\lambda$ and
$C_\lambda^t$ (resp.\ algebras $B_\lambda$ and $B_\lambda^t$)
corresponding to $J(\lambda)$ and $J^t(\lambda)$ are isomorphic via 
\begin{equation}
c \mapsto ct(\lambda),\quad c\in C_\lambda, \qquad (\mbox{ resp. } 
f(g)\mapsto f(gt(\lambda)),\quad f\in B_\lambda, g\in G).
\end{equation}
\end{remark}

The statement and proof of the next Proposition are 
analogous to (\cite{Mo}, 7).

\begin{proposition}
\label{b is semisimple}
$B_\lambda$ is a semisimple algebra (equivalently, $C_\lambda$
is a cosemisimple coalgebra).
\end{proposition}
\begin{proof}
Suppose $B_\lambda$ is not semisimple and let $I$ be the maximal non-zero
power of its radical. Then $I$ is $G$-stable and for all $z_1,z_2\in I$
and $h_1,h_2\in G$ we have $(g_1\circ z_1)(g_2\circ z_2) =0$. The last
condition is equivalent to 
\begin{equation*}
\la \Delta_\lambda(c),\,(g_1\circ z_1)\otimes (g_2\circ z_2) \ra =0
\end{equation*} 
for all $c\in C_\lambda$. In particular, we have
\begin{equation*}
\la (g_1\otimes g_2)J(\lambda),\,z_1\otimes z_2 \ra 
=\la (g_1\otimes g_2)\Delta_\lambda(1),\,z_1\otimes z_2 \ra = 0,
\end{equation*}
but since $J(\lambda)$ is invertible in $k[G]^{\otimes 2}$,
the elements $(g_1\otimes g_2)J(\lambda) (\mod A)$
span $k[G/ A]$, so that $z_1 =z_2 =0$, which is a contradiction.
\end{proof}

\begin{proposition}
\label{G acts transitively}
$G$ acts transitively on the set $\I_\lambda$ of
minimal two-sided ideals of $B_\lambda$. In particular,
all  minimal two-sided ideals of $B_\lambda$ have the same dimension.
\end{proposition}
\begin{proof}
Since $B_\lambda = F_0[G]$ as a $G$-module, the space of $G$-invariant
elements of $B_\lambda$ has dimension $1$. For any orbit of $G$ in 
$\I_\lambda$ the corresponding central idempotent of $B_\lambda$ is
$G$-invariant, so there is a single orbit.
\end{proof}

For every $\mu\in {A^*}$ consider the space $F_\mu[G]$ 
of all functions on $G$  of weight $\mu$, 
\begin{equation*}
F_\mu[G] := \{ f\in F[G] \mid f(ha) = f(h) \mu(a), \quad a\in A,\, h\in G \}.
\end{equation*}

\begin{proposition}
\label{FmuG}
$F_\mu[G]$ is a $G$-equivariant $B_{\lambda-\mu} - B_{\lambda}$ 
bimodule via 
\begin{eqnarray*}
f\circ f_\mu(g) &=& f \otimes f_\mu ((g\otimes g)J(\lambda)),
\qquad f\in B_{\lambda-\mu},\, f_\mu \in F_\mu[G], \\
f_\mu \circ f'(g) &=& f_\mu \otimes f' ((g\otimes g)J(\lambda)),
\qquad f'\in B_{\lambda},\, f_\mu \in F_\mu[G].
\end{eqnarray*}
The actions of $B_{\lambda-\mu}$ and  $B_{\lambda}^{op}$ are faithful.
\end{proposition}
\begin{proof}
First, we check that $F_\mu[G]$ is a left $B_{\lambda-\mu}$-module
and right $B_{\lambda}$-module. 

For all $f_1,f_2\in B_{\lambda-\mu}$ we have, using the definition of
$F_\mu[G]$ and dynamical twist equation~(\ref{dynamical equation}) :
\begin{eqnarray*}
f_1\circ (f_2\circ f_\mu)(g)
&=& f_1 \otimes f_2 \otimes f_\mu ((g\otimes g \otimes g)
     J^{1,23}(\lambda) J^{23}(\lambda)) \\
&=& f_1 \otimes f_2 \otimes f_\mu ((g\otimes g \otimes g) 
    J^{12,3}(\lambda) J^{12}(\lambda -\mu))\\
&=& (f_1 *_{\lambda-\mu} f_2)\circ f_\mu(g).
\end{eqnarray*}
Also, for all $f_1',f_2'\in B_{\lambda}$ we have :
\begin{eqnarray*}
(f_\mu \circ f_1')\circ f_2' (g)
&=& f_\mu \otimes f_1' \otimes f_2' ((g\otimes g \otimes g)
    J^{12,3}(\lambda) J^{12}(\lambda) ) \\
&=& f_\mu \otimes f_1' \otimes f_2' ((g\otimes g \otimes g)
    J^{1,23}(\lambda) J^{23}(\lambda)) \\     
&=& f_\mu \circ (f_1' *_{\lambda} f_2')(g).
\end{eqnarray*}
It is immediate from Equation~(\ref{twist mod A}) that the actions
of $B_{\lambda-\mu}$ and $B_{\lambda}$ commute.

To prove that the  action of  $B_{\lambda-\mu}$ is faithful
suppose that $f\circ f_\mu =0$ for some $f\in B_{\lambda-\mu}$ and
all $f_\mu \in F_\mu[G]$. Then $J(\lambda)$ being invertible
implies that $f\otimes f_\mu$ is identically equal to $0$ on 
$k[G]^{\otimes 2}$, so that $f=0$. The proof for the other action
is completely similar.
\end{proof}

\begin{proposition}
\label{all are isomorphic}
All algebras $B_\lambda$ are isomorphic to each other.
\end{proposition}
\begin{proof}
By Proposition~\ref{G acts transitively} all simple $B_\lambda$-modules 
have the same dimension, which we will denote $d_\lambda$.

Suppose $d_\lambda > d_\mu$ for some $\lambda, \mu\in A^*$. 
Let $N_\mu = \dim_k F_0[G]/d_\mu^2$ be the number of non-isomorphic
simple $B_\mu$-modules and $n_i,\,i=1,\dots,N_\mu$ be the multiplicities
with which they occur in the decomposition of $F_{\mu-\lambda}[G]$.
By Proposition~\ref{FmuG} $B_\lambda^{op}$
is isomorphic to a subalgebra of the centralizer $C_\mu$ of
$B_\mu$ in $\End_k F_{\mu-\lambda}[G]$. 
But $C_\mu = \oplus_i\, M_{n_i}(k)$, therefore $n_i \geq d_\lambda$
for all $i=1,\dots, N_\mu$. Then we have
\begin{equation*}
\dim_k F_0[G] = \sum_i\, n_i d_\mu \geq N_\mu d_\lambda d_\mu 
> N_\mu d_\mu^2 = \dim_k F_0[G],
\end{equation*}
which is a contradiction.
\end{proof}

\begin{corollary}
\label{regular}
$F_\mu[G]$ is isomorphic to the left regular $B_{\lambda-\mu}$-module
and to the right regular $B_\lambda$-module. In particular,
$B_{\lambda-\mu}$ is naturally identified with
the centralizer of $B_\lambda^{op}$ in
$\End_k(F_\mu[G])$.
\end{corollary}
\begin{proof}
In the proof of Proposition~\ref{all are isomorphic} we must have
$n_i = d_\mu$, i.e., every simple $B_\lambda$-module has the multiplicity
equal to its dimension.
\end{proof}

\begin{proposition}
\label{product of Fs}
For all $\lambda,\mu,\nu\in A^*$ the map
\begin{equation*}
\beta_{\mu\nu}^\lambda : F_\mu[G] \otimes_{B_{\lambda}}  F_\nu[G]
\to F_{\mu+\nu}[G]
\end{equation*}
defined by $\beta_{\mu\nu}^\lambda(f_\mu\otimes f_\nu)(g)
= f_\mu\otimes f_\nu ((g\otimes g)J(\lambda+\nu))$
for all $f_\mu\in F_\mu[G],\, f_\nu\in F_\nu[G],$ and $g\in G$
is an isomorphism of $G$-equivariant $B_{\lambda-\mu} - B_{\lambda+\nu}$ 
bimodules.
\end{proposition}
\begin{proof}
To show that each $\beta_{\mu\nu}^\lambda$ is well-defined,
we compute
\begin{eqnarray*}
\beta_{\mu\nu}^\lambda(f_\mu\circ f\otimes f_\nu)(g)
&=& f_\mu\circ f\otimes f_\nu((g\otimes g)J(\lambda+\nu))\\
&=& f_\mu \otimes f \otimes f_\nu((g\otimes g\otimes g)
    J^{12,3}(\lambda+\nu)J^{12}(\lambda))\\
&=& f_\mu \otimes f \otimes f_\nu((g\otimes g\otimes g)
    J^{1,23}(\lambda+\nu)J^{23}(\lambda+\nu))\\
&=& f_\mu \otimes f\circ f_\nu((g\otimes g)J(\lambda+\nu))\\
&=& \beta_{\mu\nu}^\lambda(f_\mu \otimes f\circ f_\nu)(g)
\end{eqnarray*}
for all $f\in B_\lambda$. The third equality above uses
Equation~(\ref{dynamical equation}).
It is clear that $\beta_{\mu\nu}^\lambda$ is
a $G$-module homomorphism.  To see that it is invertible we first 
observe that from Corollary~\ref{regular} we have 
$\dim_k(F_\mu[G] \otimes_{B_{\lambda}}  F_\nu[G]) = |G|/|A|$, so that
the map in question is between spaces of the same dimension. So it suffices
to show that it is surjective. 
But the range of $\beta_{\mu\nu}^\lambda$ clearly coincides with
the range of the map given by the composition of the right translation
by $J(\lambda)$ (which is surjective since $J(\lambda)$ is invertible)
and the usual multiplication map $F_\mu[G]\otimes F_\nu[G] \to F_{\mu+\nu}[G]$
(which is surjective because there are invertible functions in each
$F_\mu[G]$).
\end{proof}

\begin{definition}
\label{minimal twist}
We say that a dynamical twist $J(\lambda)$ is {\em minimal}
if $B_\lambda$ is a simple algebra for some (and hence for all) 
$\lambda\in A^*$. 
\end{definition}

\begin{remark}
\begin{enumerate}
\item[(i)]
The property of $J(\lambda)$ being minimal is invariant 
under gauge transformations.
\item[(ii)]
If $K$ is a subgroup of $G$ containing $A$ and $J(\lambda): A^*\to k[K]^2$ is
a  dynamical twist, then it is obviously a dynamical twist in $k[G]$
which is minimal if and only if  $K=G$ and $J(\lambda)$ is minimal in $k[K]$.
\end{enumerate}
\end{remark}

\begin{definition}
\label{minimizable}
A dynamical twist $J(\lambda)$ in $k[G]$ is said to be {\em minimizable} if
there exists a subgroup $K$ of $G$ containing $A$  and a minimal dynamical 
twist $J'(\lambda) : A^*\to k[K]\otimes k[K]$ such that $J(\lambda)$ is gauge
equivalent to $J'(\lambda)$, where the latter is regarded as a
dynamical twist in $k[G]$. 
\end{definition}

As it was shown in \cite{EG}, all twists are minimizable if $A =\{1\}$.
On the contrary, it is not always the case for $A \neq\{1\}$, moreover, it is
not always true that $A$ is conjugate to a subgroup of $K$,
see Example~\ref{non-minimizable data}. Also, $J(\lambda)$ is not
necessarily minimizable even if  $A=K$, see Example~\ref{affine group}.
 
\end{section}

%%%%%%%%%%%%%%%%%%%%%%%%%%%%%%%%%%%%%%%%%%%%%%%%%%%%%%%%%%%%%%%%%%%%%%%%%%%%
%%%%%%% Projective representations and dynamical data   %%%%%%%%%%%%%%%%%%%%
%%%%%%%%%%%%%%%%%%%%%%%%%%%%%%%%%%%%%%%%%%%%%%%%%%%%%%%%%%%%%%%%%%%%%%%%%%%%

\begin{section}
{Projective representations and dynamical data}

\begin{definition}
\label{projective rep}
Let $c: G\times G \to k^\times$ be a $2$-cocycle on $G$ and $V\neq 0$ 
be a vector
space over $k$. A {\em projective} representation of $G$ on $V$ 
with cocycle (Schur multiplier) $c$ is a map
$\pi : G \to GL(V)$ such that $\pi(g)\pi(h) = c(g,\,h)\pi(gh)$. A projective
representation $\pi$ is {\em linearizable} if $c$ is a coboundary and 
{\em linear} if $c\equiv 1$.  It is {\em irreducible} if $\pi(G)$ generates 
$\End_k V$ as a vector space.
\end{definition}

\begin{definition}
\label{projective isomorphism}
A {\em projective isomorphism} of two projective representations 
$\pi : G\to GL(V)$ and $\pi' : G\to GL(V')$ is an invertible
$k$-linear map  $\phi : V \to V'$ such that
\begin{equation}
\pi'(g) \circ \phi = \alpha(g) \phi \circ \pi(g), \qquad g\in G
\end{equation}
for some function $\alpha : G \to k$.
\end{definition}

It is clear that $\alpha$ takes values in $k^\times$ and that it is
uniquely defined by $\phi$. The cocycles of $\pi$ and $\pi'$ differ 
by the coboundary $\alpha(g)\alpha(h)\alpha(gh)^{-1}$.

\begin{definition}
We will call $\alpha$ the {\em multiplier} function for $\phi$. 
\end{definition}

\begin{remark}
\begin{enumerate}
\item[(a)] 
Every projective representation of $G$ descends to a homomorphism
$G\to PGL(V)$ to the projective linear group of $V$.
\item[(b)]
Let $V_i,\, i=1,2$, be projective representations of $G$.
Then $V_1\otimes V_2$ is a projective representation of $G$
in the obvious way. If the cocycles of $V_1$ and $V_2$ are the same, 
then $V_1\oplus V_2$ is also a projective representation of $G$.
\item[(c)]
For every 2-cocycle $c$ on $G$ there exists
a canonical central extension $\G$ of $G$ by a subgroup $L$ of $k^\times$
in which the cocycle $c$ takes values :
\begin{equation}
1 \longrightarrow  L \longrightarrow  \G \mathop{\longrightarrow}\limits^{p} 
G \longrightarrow  1.
\end{equation}
Explicitly, the multiplication 
in $\G$ is given by $(g,\,x)(h,\, y) = (gh,\,c(g,h)xy)$, for all $g,h\in G$
and $x,y\in L$. Any projective representation $\pi :  G \to GL(V)$ 
with Schur multiplier $c$ canonically lifts to a linear representation 
$\widehat{\pi}: \G \to GL(V)$, via
$\widehat{\pi}(g,x) = x\pi(g)$. 
In this linear representation $\widehat{\pi}$ 
the elements of $L\subset k^\times$  act on $V$ by scalar multiplication.
\item[(d)]
Projectively isomorphic representations of $G$ with the same Schur
multiplier are not necessarily isomorphic as linear representations of $\G$
even for $c=1$.  For example, any one-dimensional representation 
$\chi : G \to k^\times$ is projectively 
isomorphic to the trivial representation of $G$.
\end{enumerate}
\end{remark}

We are ready to define the main object of this paper.

\begin{definition}
\label{dynamical datum}
A {\em dynamical datum} for a pair $(G, A)$ consists of
a subgroup $K$ of $G$ and a family of projective
irreducible representations $\pi_\lambda : K\to GL(V_\lambda),\,
\lambda\in A^*$,
such that $V_{\lambda-\mu} \otimes V_\lambda^*$ is linear and
\begin{equation}
\label{product of V's}
\Ind_K^G (V_{\lambda} \otimes V_{\lambda-\mu}^*) \cong \Ind_A^G\,\mu
\end{equation}
for all $\lambda,\mu \in {A^*}$. 
%such that $\id_{V_\lambda}$ corresponds
%to $\eps\in \Ind_A^G 0 \cong F_0[G]$.
\end{definition}

\begin{remark}
Note that Definition~\ref{dynamical datum} in particular
implies that the $2$-cocycles (Schur multipliers) of all $V_\lambda$'s  
are equal, so that there is a canonical common central extension 
\begin{equation}
1\longrightarrow  L \longrightarrow \K \longrightarrow  K \longrightarrow  1,
\end{equation}
linearizing all $\pi_\lambda,\, \lambda\in A^*$. 
\end{remark}

\begin{proposition}
\label{V 's are mutually distinct}
For any dynamical datum, 
$V_\lambda$ and $V_\mu$ are non-equivalent representations of $\K$
for $\lambda\neq \mu$.
\end{proposition}
\begin{proof}
From Definition~\ref{dynamical datum} we have 
$\Hom_{\K}(V_\lambda,\,V_\mu\otimes X) = X[\lambda-\mu]$
for any $G$-module $X$ (regarded as a $\K$-module 
on which the kernel $L$ acts trivially),
where $X[\lambda]\subset X$ is a weight subspace of $X$ defined in
\ref{weight subspace}. 
Taking the trivial $G$-module $X=k$ we get $\Hom_{\K}(V_\lambda,\,V_\mu)=0$.
\end{proof}

\begin{remark}
Nevertheless, $V_\lambda$ and $V_\mu$, $\lambda\neq\mu$, can be 
projectively isomorphic as projective representations of $K$, 
e.g., when they are $1$-dimensional.
\end{remark}

\begin{definition}
\label{iso of dynamical data}
Let  $K$ and $K'$ be two subgroups of $G$.
An {\em isomorphism} of  dynamical data
$(K,\,\{V_\lambda \mid \lambda\in {A^*} \})$ and 
$(K',\,\{V_\lambda' \mid \lambda\in {A^*} \})$  
is an element $g\in G$ such that $\Ad_g\,K=K'$, where $\Ad_g$
is the adjoint action of $g$ on $G$, and a family of projective isomorphisms  
$\phi_\lambda: V_\lambda \cong V_\lambda',\,\lambda\in A^*$ 
of projective representations $\pi_\lambda$ and $\pi_\lambda'\circ \Ad_g$,
having the same multiplier function $\alpha(h)=\alpha_\lambda(h),\, h\in K$.
\end{definition}

\begin{definition}
A dynamical datum is called {\em minimal} if $K=G$.
\end{definition}

\end{section}

%%%%%%%%%%%%%%%%%%%%%%%%%%%%%%%%%%%%%%%%%%%%%%%%%%%%%%%%%%%%%%%%%%%%%%%%%%%%
%%%  A dynamical datum defined by a dynamical twist  %%%%%%%%%%%%%%%%%%%%%%%
%%%%%%%%%%%%%%%%%%%%%%%%%%%%%%%%%%%%%%%%%%%%%%%%%%%%%%%%%%%%%%%%%%%%%%%%%%%%

\begin{section}
{A dynamical datum defined by dynamical twist}
\label{4}

In this section for every dynamical twist $J(\lambda)$ we define an
isomorphism class of dynamical data, cf.\  Definitions~\ref{dynamical datum}
and \ref{iso of dynamical data}. We use the notation of Section~3.

\begin{proposition}
\label{Ilambda as Gsets}
Let $\mathcal{I}_\lambda$ be the set of minimal $2$-sided ideals
of $B_\lambda$. There are isomorphisms of $G$-sets 
$\tau_{\lambda\mu} : \mathcal{I}_\lambda \to \mathcal{I}_\mu$
such that $\tau_{\lambda\mu}\circ \tau_{\mu\nu} = \tau_{\lambda\nu}$.

In other words, we can simultaneously identify all $\mathcal{I}_\lambda$'s
as $G$-sets, $\mathcal{I}_\lambda \cong \mathcal{I}$.
\end{proposition}
\begin{proof}
Observe that a decomposition of $F_{\mu -\lambda}[G]$ into the sum 
of simple $B_\mu-B_\lambda$-bimodules establishes an isomorphisms
of $G$-sets $\tau_{\lambda\mu} : \mathcal{I}_\lambda \to \mathcal{I}_\mu$.
But the relative tensor product of bimodules $F_{\mu -\lambda}[G]$ and
$F_{\lambda-\nu}[G]$ is isomorphic to $F_{\mu -\nu}[G]$ by 
Proposition~\ref{product of Fs}, whence the composition rule of 
$\tau$'s follows.
\end{proof}

Thus, we can canonically identify stabilizers of points of 
$\mathcal{I}_\lambda$. Let us fix a point $i \in \mathcal{I}$ and let 
$K\subset G$ be the stabilizer subgroup of $i$ and $B_\lambda^i
\subset B_\lambda$ be the minimal ideal in $\mathcal{I}_\lambda$ corresponding
to $i$. Let $V_\lambda$ be a (unique up to an isomorphism) simple
$B_\lambda^i$-module, then  $K$ acts on $\End_k V_\lambda$ by automorphisms 
and we can choose  an irreducible projective representation 
\begin{equation}
\pi_\lambda: K \to GL(V_\lambda)
\end{equation}
such that the action of $g\in K$ is given by $\Ad_{\pi_\lambda(g)}$.

\begin{proposition}[cf.\ \cite{Mo}, 8]
\label{proj irreducible}
Projective representations $\pi_\lambda$ are irreducible.
\end{proposition}
\begin{proof}
We need to show that the space of $K$-invariant elements of $B_\lambda^i$ has
dimension $1$. Any such an element injectively corresponds  to a 
$G$-invariant element of $B_\lambda$, but the latter  is equivalent to 
$k[G/ A]$ as a left $G$-module, for which the space of invariants 
is $1$-dimensional.
\end{proof}

Observe that  the $K$-equivariant ${B_{\mu}^i}-{B_\nu^i}$ bimodule 
$F_{\mu\nu}^i = B_{\mu}^i F_{\nu-\mu}[G] B_\nu^i$ 
is such that $F_{\nu-\mu}[G] \cong k[G]\otimes_{k[K]} F_{\mu\nu}^i$, i.e.,
$\Ind_K^G F_{\mu\nu}^i \cong  \Ind_A^G (\mu-\nu)$.

\begin{theorem}
\label{twist to datum}
Every dynamical twist $J(\lambda): A^*\to k[G]\otimes k[G]$ defines
an isomorphism class of dynamical data. Moreover, any two
gauge equivalent dynamical twists $J(\lambda)$
and $\ti{J}(\lambda)$  define the same isomorphism class of data.
\end{theorem}
\begin{proof}
Given a dynamical twist $J(\lambda)$ fix a point $i \in \mathcal{I}$
and a subgroup $K\subset G$ as above. Consider the corresponding algebras
$B_\lambda^i$ and $B_\lambda^i-B_\mu^i$ bimodules $F_{\lambda\mu}^i$ and
fix simple $B_\lambda^i$-modules  $V_\lambda := V_\lambda^i$.  
Let $\rho_{\lambda\lambda} : F_{\lambda\lambda} \to \End_k V_\lambda$
be the isomorphisms identifying $K$-algebras and choose
bimodule isomorphisms 
$\rho_{\lambda\mu} :  F_{\lambda\mu}^i \to V_\lambda\otimes V_\mu^*$.
Let $\Sh_{\lambda\mu}(g)$ denote
the action of $g\in K$ on $F_{\lambda\mu}^i \subset F_{\mu-\lambda}[G]$
by shifts, cf.\ Equation~(\ref{action on FG}). For every $\lambda$
choose an irreducible projective representation $\pi_\lambda : K \to
GL(V_\lambda)$ such that
\begin{equation}
\label{TD 1}
\pi_\lambda(g) \otimes \pi_\lambda^*(g)^{-1} =
\rho_{\lambda\lambda}^{-1} \circ \Sh_{\lambda\lambda}(g) \circ 
\rho_{\lambda\lambda}.
\end{equation}
Then for all $x\in F_{\lambda\mu}^i,\, f\in B_\lambda^i, f'\in B_\mu^i$
we have
\begin{eqnarray*}
\lefteqn{
(\tau_{\lambda\mu} \circ \Sh_{\lambda\mu}(g)^{-1} \circ 
\rho_{\lambda\mu}^{-1}) (\pi_\lambda(g) \otimes \pi_\mu^*(g)^{-1})
(f\circ x \circ f') = }\\
&=&  f\circ (\rho_{\lambda\mu} \circ \Sh_{\lambda\mu}(g)^{-1} \circ 
\rho_{\lambda\mu}^{-1})\circ  (\pi_\lambda(g) \otimes \pi_\mu^*(g)^{-1})x 
\circ f',
\end{eqnarray*}
i.e., $(\rho_{\lambda\mu} \circ \Sh_{\lambda\mu}(g)^{-1} \circ 
\rho_{\lambda\mu}^{-1})\circ (\pi_\lambda(g) \otimes \pi_\mu^*(g)^{-1})$
commutes with the action of $B_\lambda^i \otimes (B_\mu^i)^{op}$
on $F_{\lambda\mu}^i$. Therefore, by Schur's Lemma,
\begin{equation}
\label{TD 2}
\pi_\lambda(g) \otimes \pi_\mu^*(g)^{-1} =\gamma_{\lambda\mu}(g)
\rho_{\lambda\mu}^{-1} \circ \Sh_{\lambda\mu}(g) \circ \rho_{\lambda\mu}
\end{equation}
for some functions $\gamma_{\lambda\mu} : K \to k^\times$.
By Proposition~\ref{product of Fs}, 
$F_{\lambda\mu}^i \otimes_{B_\mu^i} F_{\mu\nu}^i \cong F_{\lambda\nu}^i$
as $K$-equivariant bimodules, whence $\gamma_{\lambda\mu}(g)
\gamma_{\mu\nu}(g) = \gamma_{\lambda\nu}(g)$ for all $\lambda,\mu,\nu$.
Replacing $\pi_\lambda(g)$ by a projectively isomorphic representation
$\pi_\lambda'(g) = a_\lambda(g)\pi_\lambda(g)$, where 
$a_\lambda(g)^{-1} a_\mu(g) = \gamma_{\lambda\mu}(g)$,
we obtain a collection of irreducible projective representations $\pi_\lambda'$
of $K$ such that $V_\lambda \otimes V_\mu^*$ is linear and induces
$F_{\mu-\lambda}[G] \cong \Ind_A^G(\lambda-\mu)$, i.e., a dynamical datum
for $(G, A)$.

Let us show that gauge equivalent twists define isomorphic dynamical data,
for any choice of projective representations $\pi_\lambda$. This will
imply, in particular, that the isomorphism class of a dynamical datum
we have constructed above is well defined.

Let $\ti{J}(\lambda)$ be a dynamical twist in $k[G]$ gauge equivalent
to $J(\lambda)$ via a gauge transformation $t(\lambda)$. 
By Remark~\ref{ge --> iso of Bl} the map
$t_\lambda : f(g)\mapsto f(gt(\lambda))$ defines an isomorphism
between the corresponding $G$-algebras $B_\lambda$ and $\ti{B}_\lambda$.
This establishes a bijective correspondence between minimal ideals
of $B_\lambda^i$ and $\ti{B}_\lambda^i$ of these algebras and allows
to identify their stabilizers, so we can assume that both
$J(\lambda)$ and $\ti{J}(\lambda)$ define the same subgroup $K$.

Let $\ti{F}_{\lambda\mu}^i$ be the corresponding 
$\ti{B}_\lambda^i-\ti{B}_\mu^i$ bimodules. Fix vector spaces
$\ti{V}_\lambda:=\ti{V}_\lambda^i$ and  identify $K$-algebras
$\ti{B}_\lambda^i$ and $\End_k \ti{V}_\lambda$. Choose isomorphisms
$\ti{\rho}_{\lambda\mu} :  \ti{F}_{\lambda\mu}^i \to \ti{V}_\lambda\otimes 
\ti{V}_\mu^*$ of bimodules, and projective representations
$\ti{\pi}_\lambda: K \to GL(\ti{V}_\lambda)$ such that
\begin{equation*}
\ti{\pi}_\lambda(g) \otimes \ti{\pi}_\lambda^*(g)^{-1} =
\ti{\rho}_{\lambda\lambda}^{-1} \circ \ti{\Sh}_{\lambda\lambda}(g) \circ 
\ti{\rho}_{\lambda\lambda},
\end{equation*}
where $\ti{\Sh}_{\lambda\lambda}(g)$ denotes the action of $g$ on
$\ti{F}_{\lambda\mu}^i$. Arguing as above, we get functions
$\ti{\gamma}_{\lambda\mu}(g)$ and $\ti{a}_\lambda(g)$ such that
\begin{equation}
\label{TD 3}
\ti{\pi}_\lambda(g) \otimes \ti{\pi}_\mu^*(g)^{-1} =\ti{\gamma}_{\lambda\mu}(g)
\ti{\rho}_{\lambda\mu}^{-1} \circ \ti{\Sh}_{\lambda\mu}(g) \circ 
\ti{\rho}_{\lambda\mu},
\end{equation}
and $\ti{\pi}_\lambda'(g) = \ti{a}_\lambda(g)\ti{\pi}_\lambda(g)$
(where $\ti{a}_\lambda(g)  \ti{a}_\mu(g)^{-1} = 
\ti{\gamma}_{\lambda\mu}(g)^{-1}$)
make $\ti{V}_\lambda\otimes \ti{V}_\mu^*$ linear.

Isomorphisms between simple $K$-algebras $B_\lambda^i$ and $\ti{B}_\lambda^i$
yield projective isomorphisms $\phi_\lambda$ (defined up to scaling)
between $\pi_\lambda$ and $\ti{\pi}_\lambda$, $\lambda\in A^*$:
\begin{equation}
\label{TD 4}
\ti{\pi}_\lambda(g) \circ \phi_\lambda =
\alpha_\lambda(g) \phi_\lambda \circ \pi_\lambda(g), \qquad g\in K,
\end{equation}
for some multipliers $\alpha_\lambda: K \to k^{\times}$.

Comparing Equations~(\ref{TD 2}) and (\ref{TD 3}) and using
(\ref{TD 4}) we conclude that
\begin{equation*}
\frac{\ti{\gamma}_{\lambda\mu}(g)}{\gamma_{\lambda\mu}(g)} =
\frac{\alpha_\lambda(g)}{\alpha_\mu(g)}, \qquad g\in K.
\end{equation*}
We have $\ti{\pi}_\lambda'(g) \circ \phi_\lambda =
\alpha_\lambda'(g) \phi_\lambda \circ \pi_\lambda'(g)$, where
$\alpha_\lambda'(g) = \alpha_\lambda(g) \frac{\ti{a}_\lambda(g)}{a_\lambda(g)}$
and
\begin{equation*}
\frac{\alpha_\lambda'(g)}{\alpha_\mu'(g)} =
\frac{\alpha_\lambda(g)}{\alpha_\mu(g)}
\frac{\ti{a}_\lambda(g)}{\ti{a}_\mu(g)}
\frac{a_\mu(g)}{a_\lambda(g)} = 1,
\end{equation*}
so that all multipliers $\alpha_\lambda'$ are equal, i.e., dynamical data
constructed from $J(\lambda)$  and $\ti{J}(\lambda)$ are isomorphic.
\end{proof}

\begin{corollary}
The above construction assigns minimal dynamical data to minimal
dynamical twists.
\end{corollary}
\begin{proof}
This is clear since $B_\lambda$ is simple if and only if $G=K$.
\end{proof}

\begin{remark}
\label{composition1}
It follows from the definition of a dynamical datum that for every
$\mu\in A^*$ the weight subspace $F_\mu[G]$ can be regarded as the space 
of functions on $G/K$ with values in $\Hom_k(V_\lambda,V_{\lambda-\mu})$,
or, equivalently, as the space of $K$-homomorphisms from $k[G]$ to
$\Hom_k(V_\lambda,V_{\lambda-\mu})$, where $h\in K$ acts on $k[G]$
by the right multiplication by $h^{-1}$ and by conjugation on
$\Hom_k(V_\lambda,V_{\lambda-\mu})$.

With this identification the map $F_\nu[G]\otimes_{B_{\lambda-\mu}} F_\mu[G]
\to F_{\nu+\mu}[G]$ defined in Proposition~\ref{product of Fs} 
is given by the composition
\begin{equation}
\label{star}
\begin{split}
k[G]\to k[G]\otimes k[G] 
& \to  \Hom_k(V_{\lambda-\mu},V_{\lambda-\mu-\nu})  
\otimes \Hom_k(V_\lambda,V_{\lambda-\mu})  \to  \\
& \to \Hom_k(V_\lambda,V_{\lambda-\mu-\nu}),
\end{split}
\end{equation}
where the first map is the comultiplication $g\to g\otimes g$, the second uses
the above description of the weight subspaces, and the third is the composition
of homomorphisms.
\end{remark}

\end{section}

%%%%%%%%%%%%%%%%%%%%%%%%%%%%%%%%%%%%%%%%%%%%%%%%%%%%%%%%%%%%%%%%%%%%%%%%%%%%% 
%%%  Construction of a dynamical twist from dynamical datum  %%%%%%%%%%%%%%%% 
%%%%%%%%%%%%%%%%%%%%%%%%%%%%%%%%%%%%%%%%%%%%%%%%%%%%%%%%%%%%%%%%%%%%%%%%%%%%%

\begin{section}
{Construction of a dynamical twist from dynamical datum}
\label{5}

We use the exchange construction that appeared 
in \cite{EV} (see also \cite{ES}) to produce a dynamical twist from 
a dynamical datum for $(G, A)$. 

Given a dynamical datum $(K,\,\{V_\lambda \mid \lambda\in {A^*} \})$
as in Definition~\ref{dynamical datum}
let $\K$ be the canonical central extension of $K$ linearizing all $V_\lambda$.
Then every $G$-module $X$ is also a $\K$-module on which the kernel of
the extension (which is a central subgroup of $\K$) acts trivially.

Choose isomorphisms $\epsilon_{\lambda\mu}:
V_\lambda \otimes V_{\lambda-\mu}^* \cong
\Ind_A^G\mu$ such that for $\lambda=\mu$ the identity of
$\End_k V_\lambda$ is mapped to $1$. 
Using the definition of dynamical datum and
Frobenius reciprocity for induced modules we have
\begin{equation}
\label{eq 20}
\begin{split}
\Hom_{\K}(V_\lambda,\, V_{\lambda-\mu}\otimes X) &\cong 
\Hom_{K}(V_\lambda \otimes V_{\lambda-\mu}^*,\, X) \\
&\cong \Hom_G(\Ind_A^G\mu, X) \cong X[\mu],
\end{split}
\end{equation}
where the second isomorphism is defined by $\epsilon_{\lambda\mu}$
and the other two isomorphisms are canonical.

Let us denote $\Psi(\lambda, x)$
the homomorphism $V_\lambda \to V_{\lambda-\mu}\otimes X$
determined by $x\in X[\mu]$ through the above sequence of
isomorphisms.
Let $Y$ be another $G$-module and $y\in Y[\nu],\,\nu\in A^*$. Then
\begin{equation}
(\Psi(\lambda-\mu,y)\otimes \id)\circ \Psi(\lambda, x) :
V_\lambda \to V_{\lambda-\mu-\nu}\otimes Y\otimes X
\end{equation}
corresponds to a unique element in $(Y\otimes X)[\mu+\nu]$.
We introduce a linear operator $J_{Y,X}(\lambda):Y\otimes X \to Y\otimes X$
by setting
\begin{equation}
(\Psi(\lambda-\mu,y)\otimes \id)\circ \Psi(\lambda, x):=
\Psi(\lambda, J_{Y,X}(\lambda)(y\otimes x)).
\end{equation}
We will show that the function $J(\lambda)\in k[G]\otimes k[G]$
defined by this equation is a dynamical twist for $k[G]$.

\begin{remark}
\label{composition2}
Note that the above defined $J(\lambda)$ is related to the composition map
\begin{equation*}
\Hom_k(V_{\lambda-\mu},V_{\lambda-\mu-\nu}) \otimes
\Hom_k(V_\lambda,V_{\lambda-\mu})
 \to \Hom_k(V_\lambda,V_{\lambda-\mu-\nu}),
\end{equation*}
in the following way. Take $f_\mu\in F_\mu[G],\, f_\nu\in F_\nu[G]$.
If we regard $\Psi(\lambda, f_\mu)$ as a $K$-homomorphism from $k[G]$
to $\Hom_k(V_\lambda,V_{\lambda-\mu})$ (where $h\in K$ acts on $k[G]$
as $h\cdot g = gh^{-1}$) then
\begin{equation}
\label{comp}
\Psi(\lambda, J(\lambda)(f_\nu\otimes f_\mu)) =
\comp(\Psi(\lambda-\mu,f_\nu)\otimes \Psi(\lambda, f_\mu))\circ\Delta,
\end{equation}
where $\Delta: k[G]\to k[G]\otimes k[G]$ is the comultiplication and
$\comp()$ is the composition of homomorphisms (cf.\  Equation~(\ref{star})).
\end{remark}

\begin{lemma}
$J(\lambda)$ has zero weight.
\end{lemma}
\begin{proof}
For all $G$-modules $X,Y$ and elements
$x\in X[\mu]$, $y\in Y[\nu]$, and $a\in A$ we have :
\begin{equation*}
J_{Y,X}(\lambda)(ay\otimes ax) = (\mu+\nu)(a) J_{Y,X}(\lambda)(y\otimes x)
= (a\otimes a) J_{Y,X}(\lambda)(y\otimes x),
\end{equation*}
since $J_{Y,X}(\lambda)$ maps $Y[\nu] \otimes X[\mu]$ to
$(Y\otimes X)[\mu+\nu]$.
\end{proof}

\begin{lemma}[cf.\ \cite{Mo}]
$J(\lambda)$ is invertible for all $\lambda\in A^*$.
\end{lemma}
\begin{proof}
We need to show that 
\begin{equation*}
J(\lambda) : \bigoplus_\mu Y[\mu-\eta] \otimes X[\lambda-\mu]
\to  (Y\otimes X)[\lambda-\eta]
\end{equation*}
is surjective for all $\lambda,\eta\in A^*$, where $X,Y$ are
$G$-modules. Equivalently, we need to prove that the composition map
\begin{equation*}
\bigoplus_\mu \Hom_{\K}(V_\lambda\otimes X^*,\,V_\mu) \otimes
\Hom_{\K}(V_\mu, V_\eta\otimes Y)
\to  \Hom_{\K}( V_\lambda\otimes X^*,\, V_\eta\otimes Y)
\end{equation*}
is surjective. But this follows from the fact that $V_\eta\otimes Y
\cong \oplus_\mu Y[\mu-\eta]\otimes V_\mu$ as $\K$-modules.
Indeed, $\Hom_{\K}(V_\mu,\,V_\eta\otimes Y)\cong Y[\mu-\eta]$ by 
Equation~(\ref{eq 20}) and $V_\mu's$ are mutually non-equivalent by 
Proposition~\ref{V 's are mutually distinct}. Therefore, a copy of
$\oplus_\mu Y[\mu-\eta]\otimes V_\mu$ is contained in $V_\eta\otimes Y$,
and has the same dimension.
\end{proof}

\begin{theorem}
\label{J is a twist}
$J(\lambda)$ is a dynamical twist for $k[G]$.
Moreover, two isomorphic dynamical data
$(K,\,\{V_\lambda \mid \lambda\in {A^*} \})$ and 
$(K',\,\{V_\lambda' \mid \lambda\in {A^*} \})$ define
the same gauge equivalence class of dynamical twists.
\end{theorem}
\begin{proof}
We have already seen that $J(\lambda)$ is an invertible zero weight
function, we need to show that it satisfies the twist properties of
Definition~\ref{dynamical twist}. Let $X,Y,Z$ be $G$-modules and 
$x\in X[\mu],\, y\in Y[\nu],\, z\in Z[\eta]$. Consider the composition
of $\K$-module homomorphisms
\begin{equation*}
V_\lambda \to V_{\lambda-\mu} \otimes X \to 
V_{\lambda-\mu-\nu}\otimes Y \otimes X \to 
 V_{\lambda-\mu-\nu-\eta} \otimes Z\otimes Y \otimes X
\end{equation*}
defined by $x,y,z$. The associativity law gives two different ways
of factorizing it:
\begin{eqnarray*}
\lefteqn{(\Psi(\lambda-\mu,J_{Z,Y}(\lambda-\mu)(z\otimes y))\otimes\id)
\circ \Psi(\lambda,x)  = }\\
&= & (\Psi(\lambda-\mu, z)\otimes\id)\circ 
\Psi(\lambda, J_{Y,X}(\lambda)(y\otimes x)),
\end{eqnarray*}
whence we have $J_{Z\otimes Y, X}(\lambda)(J_{Z,Y}(\lambda-\mu)\otimes\id)
= J_{Z,Y\otimes X}(\lambda) (\id\otimes J_{Y,X}(\lambda))$, i.e.,
$$
J^{12,3}(\lambda) J^{12}(\lambda-h^{(3)}) =
J^{1,23}(\lambda) J^{23}(\lambda)
$$ 
for all $\lambda$.
The counital properties  of $J(\lambda)$ are clear since 
$\Psi(\lambda, \eps) = \id_{V_\lambda}$.

In proving the gauge equivalence of dynamical twists coming from isomorphic 
data we  may assume that $K=K'$.
Let $\phi_\lambda : V_\lambda \to V_\lambda',\, \lambda\in A^*$ 
be an  isomorphism between two dynamical data, $\Psi(\lambda, x),\, 
\Psi'(\lambda,x)\,(x\in X[\mu])$ be associated homomorphisms of modules,
and $J(\lambda),\, J'(\lambda)$ the dynamical twists.
Define $t_X(\lambda) : X\to X$ by the identity
\begin{equation*}
\Psi'(\lambda, t_X(\lambda)x) = (\phi_{\lambda-\mu}\otimes \id) \circ
\Psi(\lambda, x) \circ \phi_{\lambda}^{-1}.
\end{equation*}
Note that $t_X(\lambda)$ is well defined since all $\phi_\lambda$
have the same multiplier function, cf.\ Definition~\ref{iso of dynamical data}.
Then $t_X(\lambda)$ is a zero weight function taking invertible values
and satisfying $\eps(t_X(\lambda))=1$ since $t_k(\lambda) =\id_k$.
For any  $y\in Y[\nu]$ we compute
\begin{eqnarray*} 
\lefteqn{\Psi'(\lambda, J_{Y,X}'(\lambda)(y\otimes x)) = }\\
&=& (\Psi'(\lambda-\mu,y)\otimes \id)\circ \Psi'(\lambda,x)\\
&=& (\phi_{\lambda-\mu-\nu} \otimes \id)\circ 
    \Psi(\lambda-\mu,t_Y(\lambda-\mu)y) \circ
    \Psi(\lambda,t_X(\lambda)(x))\circ \phi_{\lambda}^{-1}\\
&=& (\phi_{\lambda-\mu-\nu} \otimes \id)\circ
    \Psi(\lambda,J_{Y,X}(\lambda)(t_Y(\lambda-\mu)y\otimes t_X(\lambda)x)
    \circ \phi_{\lambda}^{-1}\\
&=& \Psi'(\lambda, t_{Y\otimes X}(\lambda)^{-1} 
     J_{Y,X}(\lambda)(t_Y(\lambda-\mu)y\otimes t_X(\lambda)x)),
\end{eqnarray*}
therefore $J_{Y,X}'(\lambda) =  t_{Y\otimes X}(\lambda)^{-1} J_{Y,X}(\lambda)
(t_Y(\lambda-\mu) \otimes t_X(\lambda))$, and
\begin{equation*}
J'(\lambda) = \Delta(t(\lambda))^{-1} J(\lambda) 
(t(\lambda-h^{(2)})\otimes t(\lambda)),
\end{equation*}
i.e., $J'(\lambda)$ is a gauge transformation of $J(\lambda)$.
\end{proof}

\begin{remark}
One cannot canonically construct a concrete $J(\lambda)$ from 
$(K,\,\{ V_\lambda\})$
or vice versa. It is only possible to assign an equivalence
class of dynamical twists to every $(K,\,\{ V_\lambda\})$ and
isomorphism class of data to every $J(\lambda)$.
\end{remark}

\begin{theorem}
\label{our lovely correspondence}
The maps $D$ and $T$ between gauge equivalence classes of twists
and isomorphism classes data  described by constructions of Section~\ref{4}
and the beginning of this Section, respectively, 
\begin{equation}
\left\{
\begin{array}{c}
\mbox{gauge equivalence}\\
\mbox{classes of} \\
\mbox{dynamical twists} \\
J(\lambda): A^*\to k[G]\otimes k[G]
\end{array}
\right\}
\begin{array}{c}
{}\\
\stackrel{D}{\longrightarrow}\\
\stackrel{T}{\longleftarrow}\\
{}
\end{array}
\left\{
\begin{array}{c}
\mbox{isomorphism}\\
\mbox{classes of} \\
\mbox{dynamical data}\\
 (K,\,\{ V_\lambda\mid \lambda\in A^* \})
\end{array}
\right\},
\end{equation}
are inverses of each other, i.e., define a
bijective correspondence between the two sets involved.
\end{theorem}
\begin{proof}
Let $J(\lambda)$ be a dynamical twist, $(K,\,\{ V_\lambda\})$ be 
(a representative of the class of)
dynamical data associated to  $J(\lambda)$, and $\tilde{J}(\lambda)$
be a dynamical twist coming from the exchange construction
for $(K,\,\{V_\lambda\})$. According to Remarks~\ref{composition1} and
\ref{composition2} (cf. Equations~(\ref{star}) and (\ref{comp}))
both $J(\lambda)$ and $\tilde{J}(\lambda)$
are determined by the map
\begin{equation*}
F_\nu[G]\otimes_{B_{\lambda-\nu}} F_\mu[G] \to F_{\nu+\mu}[G]
\end{equation*}
corresponding to the composition
\begin{equation*}
\Hom_k(V_{\lambda-\mu},V_{\lambda-\mu-\nu})
 \otimes \Hom_k(V_\lambda,V_{\lambda-\mu}) 
 \to \Hom_k(V_\lambda,V_{\lambda-\mu-\nu}).
\end{equation*}
Comparing two $G$-module isomorphisms 
$\Hom_K(k[G],\Hom_k(V_\lambda,V_{\lambda-\mu})) 
\cong F_\mu[G]$ (the one established in Theorem~\ref{twist to datum}, 
corresponding to the construction of $\{ V_\lambda\}$ from $J(\lambda)$,
and the one used in the exchange construction above)
we get a zero weight function  
with invertible values $t(\lambda): A^*\to k[G]$
that implements the corresponding $G$-module automorphism
\begin{equation*}
F_\mu[G] \to F_\mu[G] : f(g) \mapsto f(gt(\lambda)), \qquad g\in G,
\end{equation*}
for all $\lambda$, and is such that
\begin{equation*}
f_\nu\otimes f_\mu((gt(\lambda-\mu)\otimes gt(\lambda))\tilde{J}(\lambda))
= f_\nu\otimes f_\mu((g\otimes g) J(\lambda)\Delta(t(\lambda))),
\end{equation*}
for all $f_\nu\in F_\nu[G],\, f_\mu\in F_\mu[G]$. Therefore
$J(\lambda)$ and $\tilde{J}(\lambda)$ are gauge equivalent,
i.e., $T\circ D =\id$, in particular $T$ is surjective.

Let us show that $T$ is also injective. Let $(K,\,\pi_\lambda : K \to
GL(V_\lambda))$ and  $(\ti{K},\,\ti{\pi}_\lambda : \ti{K} \to
GL(\ti{V}_\lambda))$ be two sets of dynamical data that produce 
dynamical twists $J(\lambda)$ and $\ti{J}(\lambda)$ gauge equivalent
to each other. By Remark~\ref{ge --> iso of Bl} their
$G$-algebras $B_\lambda$ and $\ti{B}_\lambda$ are isomorphic via
\begin{equation*}
f(g) \mapsto f(gt(\lambda)),\qquad f\in B_\lambda = F_0[G].
\end{equation*}
Also, for the corresponding bimodules $F_{\lambda\mu}$ and
$\ti{F}_{\lambda\mu}$ it follows from the exchange construction
of a dynamical twist in the beginning of this section that
there are canonical isomorphisms of $G$-modules 
$$
F_{\lambda\mu} \cong 
(\Hom_k(V_\lambda, V_{\lambda-\mu}) \otimes F[G])^K, \qquad
(\mbox{resp. }
\ti{F}_{\lambda\mu} \cong (\Hom_k(\ti{V}_\lambda, \ti{V}_{\lambda-\mu}) 
\otimes F[G])^{\ti{K}}),
$$
where $K$ (resp.\ $\ti{K}$) acts by right translations on $F[G]$ and
in a standard way on the $\Hom$ space, whereas $g\in G$ acts on $F[G]$ by the 
left translation by $g^{-1}$ and trivially on the first factor.

For $\lambda=\mu$ this implies that all the matrix
blocks of $B_\lambda$  (resp.\ $\ti{B}_\lambda$) are canonicaly isomorphic 
to  $\End_k V_\lambda$ (resp. $\End_k \ti{V}_\lambda$) and form
a $G$-homogeneous space isomorphic to $G/K$ (resp.\  $G/\ti{K}$).
Therefore $\ti{K}$ is conjugate to $K$ and we can assume $\ti{K} = K$.

Thus, we have a canonical $K$-algebra isomorphism $t_\lambda:
\End_k V_\lambda \to \End_k \ti{V}_\lambda$ and so there is a
projective isomorphism $\phi_\lambda : V_\lambda \to
\ti{V}_\lambda$ such that $\Ad\phi_\lambda =t_\lambda$.
Note that $\phi_\lambda$ is defined up to a scalar, therefore
its multiplier function $\alpha_\lambda(g)$,
\begin{equation*}
\ti{\pi}_\lambda(g) \circ \phi_\lambda =  
\alpha_\lambda(g) \phi_\lambda  \circ  \pi_\lambda(g),
\end{equation*}
is uniquely defined. 

The gauge transformation $t(\lambda)$ also defines
isomorphisms of simple $K$-equivariant bimodules :
\begin{equation*}
t_{\lambda\mu} : F_{\lambda\mu}\to \ti{F}_{\lambda\mu} :
f(g) \mapsto f(gt(\lambda)),
\end{equation*}
where $ F_{\lambda\mu}, \ti{F}_{\lambda\mu} \subset F_{\mu\lambda}[G]$,
in particular, $t_{\lambda\lambda} =t_\lambda$. As it was observed
above, isomorphisms $\rho_{\lambda\mu}: F_{\lambda\mu}\to
V_\lambda \otimes V_\mu^*$ and $\ti{\rho}_{\lambda\mu}: \ti{F}_{\lambda\mu}\to
\ti{V}_\lambda \otimes \ti{V}_\mu^*$ of equivariant $K$-bimodules
are also canonical. By Schur's Lemma, for fixed $\lambda,\mu$ we have
\begin{equation*}
( \phi_\lambda \otimes (\phi_\mu^*)^{-1}) =
C\, \ti{\rho}_{\lambda\mu} \circ t_{\lambda\mu} \circ \rho_{\lambda\mu}^{-1},
\end{equation*}
for some constant $C$,
whence replacing $\phi_\lambda$ by $C\phi_\lambda$ we get $\alpha_\lambda(g)
= \alpha_\mu(g)$, i.e., the system $\{ \phi_\lambda \}$ gives an isomorphism 
between the dynamical data in question.
\end{proof}

\begin{remark}
In the case $A=\{1\}$ Theorem~\ref{our lovely correspondence} recovers
the result of \cite{Mo}, \cite{EG}, where usual twists in $k[G]$ 
(modulo gauge equivalence) were shown
to be in  bijection with single irreducible projective representations of
subgroups of $G$ (modulo conjugation).
\end{remark}

\begin{corollary}
\label{J is minimal}
A dynamical datum is minimal, i.e., $G=K$, if and only if
the corresponding dynamical twist $J(\lambda)$ is minimal.
\end{corollary}
%
%\begin{proof}
%The second part of the proof of Theorem~\ref{our lovely correspondence}
%shows that in this case the stabilizer of any minimal ideal of $B_\lambda$
%is the group $G$ itself, i.e., $B_\lambda$ is simple.
%\end{proof}

\begin{example}
\label{non-minimizable data}
Here is an example of a dynamical datum that gives rise to a non-minimizable
dynamical twist.

Let $f :A_2 \to A_1$ be an isomorphism between two abelian subgroups 
of $G$ which are {\em not} conjugate to one another and such that for 
any irreducible character $\lambda$ of $A_1$ we have 
$\Ind_{A_1}^G\, \lambda \cong
\Ind_{A_2}^G\, f^*\lambda$. Then
taking $A=A_1$, $K=A_2$, $V_\lambda =k_\lambda$, and 
$\pi_\lambda =  f^*\lambda,\, \lambda\in A^*$,
we get a dynamical datum $(K,\,\{V_\lambda\})$ that gives rise to a 
non-minimizable dynamical twist 
(since $A$ is not contained in any subgroup conjugate to $K$).

Below is an example of such a situation that we learned
from R.~Guralnick. Take $G=S_6$, the symmetric group of degree $6$,
and let $A_1$ and $A_2$ be two non-cyclic subgroups of order $4$ with
$A_1$ moving precisely $4$ points (in a single orbit) and $A_2$
having three orbits of size $2$, e.g., 
$A_1 =\{ 1,\, (12)(34),\,(13)(24),\, (14)(23)\}$ and
$A_2 =\{ 1,\, (12)(34),\,(34)(56),\, (12)(56)\}$.
Clearly, $A_1$ and $A_2$ are not conjugate as they have orbits of different 
size. 

On the other hand, we have 
$\Ind_{A_1}^G\,1 \cong \Ind_{A_2}^G\, 1$ since $A_1$ and $A_2$
intersect each conjugacy class in $G$ in the same number of elements
(namely, the identity and $3$ conjugates of $(12)(34)$). Let $\lambda_i$
denote $3$ non-trivial irreducible characters of $A_1$. Note that
$\Ind_{A_1}^G\, \lambda_i \cong \Ind_{A_1}^G\, \lambda_j$ since the
three proper subgroups of $A_1$  (and $A_2$) are conjugate in $G$. Next,
\begin{equation*}
\Ind_{\{1\}}^G\, 1 = \Ind_{A_1}^G (\Ind_{\{1\}}^{A_1}\,1)
= \Ind_{A_1}^G\, 1 + \sum \Ind_{A_1}^G\, \lambda_i
\cong \Ind_{A_1}^G\, 1 + 3\, \Ind_{A_1}^G\, \lambda,
\end{equation*}
where $\lambda$ is any of $\lambda_i$. We have precisely the same equation
for $A_2$, and since $\Ind_{A_1}^G\,1 \cong \Ind_{A_2}^G\, 1$ as noted above,
this implies that $\Ind_{A_1}^G\, \lambda \cong \Ind_{A_2}^G\, f^*\lambda$.
\end{example}

\begin{example}
\label{affine group}
Here we describe a general construction of dynamical data with $K=A$
and compute corresponding dynamical twists in an important concrete situation.
Let $A$ be an abelian subgroup of $G$ and $f: A^* \to A^*$ be a bijection
(which is not necessarily a group isomorphism). Suppose that for all
$\lambda, \mu \in A^*$ the character $f(\lambda)-f(\mu)$ is conjugate
to  $\lambda-\mu$ via some element $g(\lambda,\mu)$ 
in the normalizer $N(A)$ of $A$, i.e.,
\begin{equation*}
(f(\lambda)-f(\mu))(a) = (\lambda-\mu)(\Ad_{g(\lambda,\mu)}a), \qquad a\in A,
\end{equation*}
where $\Ad$ denotes the left adjoint action of $G$ on itself.

Take $K=A$, $V_\lambda =k$, and $\pi_\lambda = f(\lambda),\, \lambda\in A^*$.
We claim that $D_f := (K,\, \{V_\lambda\})$ is a dynamical datum for $(G, A)$.
Indeed, we have $\Ind_A^G (\lambda-\mu)  \cong \Ind_A^G(f(\lambda)-f(\mu))$
since the characters  $\lambda-\mu$ and $f(\lambda)-f(\mu)$ are conjugate 
by $g(\lambda,\mu)\in G$.

Moreover, for any two bijections $f_1,\,f_2: A^*\to K^*$ with the above 
property,
the corresponding data $D_{f_1}$ and   $D_{f_2}$ are isomorphic if and only
if $f_1(\lambda) = \Ad_g\,f_2(\lambda) +\mbox{const}$, where $g\in N(A)$. 
In particular, $D_f$ is
isomorphic to the trivial datum (which corresponds to the constant twist
$1\otimes 1$)  if and only if  it defines a minimizable twist if and only if
$f(\lambda) = \Ad_g\,\lambda + \mbox{const},\, g\in N(A)$.

Choose elements  $g(\lambda,\mu)\in G$ and
let $X$ be a $G$-module and $x\in X[\mu]$. 
Define a corresponding $K$-module homomorphism by
\begin{equation}
\Psi(\lambda, x)\, :\, V_\lambda \to V_{\lambda-\mu} \otimes X \,:\,
1\mapsto 1 \otimes g(\lambda,\lambda-\mu)^{-1}x
\end{equation}
(note that $g(\lambda,\lambda-\mu)^{-1}$ maps $X[\mu]$ to 
$X[f(\lambda)-f(\lambda-\mu)]$ for all $\lambda,\mu\in A^*$).

Applying the exchange construction we obtain a dynamical twist
\begin{equation}
\label{vot on twist}
\begin{split}
J(\lambda) =\sum_{\mu\nu\in A^*}\, &
g(\lambda, \lambda-\mu-\nu) g(\lambda-\mu, \lambda-\mu-\nu)^{-1}  P_\nu
\otimes  \\
& \quad \otimes
g(\lambda, \lambda-\mu-\nu) g(\lambda, \lambda-\mu)^{-1}    P_\mu,
\end{split}
\end{equation}
where $P_\mu = \frac{1}{|A|}\sum_{a\in A^*}\, \mu(a)a$ 
is the projection on $k[G][\mu]$.

As a  concrete example of such a situation consider the group
$G = GL_n(F_p) \ltimes F_p^n$, where $F_p$ is the field of $p$ elements
and $A=F_p^n$ is regarded as an additive group on which $GL_n(F_p)$ acts
by multiplication, and  an arbitrary bijection $f : A^*\to A^*$ 
(note that any two non-trivial characters of $A$ are conjugate 
since $GL_n(F_p)$ acts  on $(F_p^n)^{*}\backslash\{0\}$ transitively).

We compute dynamical twist (\ref{vot on twist}) when $n=1$. In this case
$G = F_p^\times \ltimes F_p$, where $F_p^\times$ acts on $A=F_p$ by
multiplication (this $G$ can be thought as the group of affine
transformations of the line over $F_p$). Let $\zeta$ be a primitive
$p$th root of unity in $k$ and $\tau: F_p^*\to F_p$ be
the isomorphism given by
\begin{equation*}
\lambda(a) = \zeta^{\tau(\lambda)a}, \qquad
\lambda\in A^*,\,a\in A.
\end{equation*}
In what follows we suppress $\tau$ and identify $F_p^*$ and $F_p$.
For any bijection $f: F_p \to F_p$  define elements of $F_p^\times$ : 
\begin{equation}
g(\lambda,\mu) =
\begin{cases} 
\frac{f(\lambda)-f(\mu)}{\lambda-\mu}, &\mbox{if } \lambda\neq \mu;\\ 
1, &\mbox{otherwise},
\end{cases}
\end{equation}
conjugating  $\lambda-\mu$ and $f(\lambda)-f(\mu)$ for all $\lambda,\mu$.
Taking these $g(\lambda,\mu)$ and identifying $F_p^*\backslash\{0\}$ 
and  $F_p\backslash\{0\}$ we get an explicit formula 
for a dynamical twist:
\begin{equation}
\begin{split}
J(\lambda) =   \frac{1}{p^2}
\sum_{\mu\nu\in F_p^*, \, ab\in F_p}\, 
\zeta^{\nu a+\mu b}\,
& \frac{\nu}{f(\lambda-\mu)- f(\lambda-\mu-\nu)}\,
\frac{f(\lambda)-f(\lambda-\mu-\nu)}{\mu+\nu}\, a
\otimes \\
& \otimes
\frac{\mu}{f(\lambda)- f(\lambda-\mu)}\,
\frac{f(\lambda)-f(\lambda-\mu-\nu)}{\mu+\nu}\, b,
\end{split}
\end{equation}
where the fractions with zero denominators are replaced by $1$.
\end{example}

\end{section} 

%%%%%%%%%%%%%%%%%%%%%%%%%%%%%%%%%%%%%%%%%%%%%%%%%%%%%%%%%%%%%%%%%%%%%%%%%%%%
%%%%%%%%%%%%  A method for finite nilpotent groups %%%%%%%%%%%%%%%%%%%%%%%%%
%%%%%%%%%%%%%%%%%%%%%%%%%%%%%%%%%%%%%%%%%%%%%%%%%%%%%%%%%%%%%%%%%%%%%%%%%%%%

\begin{section}
{A method of constructing dynamical twists for finite nilpotent groups}

Here we describe how to produce dynamical data from nilpotent
Lie algebras over finite fields.
Let $\g^0$ be a split semisimple Lie algebra over the field $F=F_p$ of
$p$ elements (e.g., $\mathfrak{sl}_n(F_p)$). 
Consider a nilpotent Lie algebra 
$\g = \g^0tF[t]/t^{n+1}$ defined by $[xt^i, yt^j] = [x,y] t^{i+j}$
for all $x,y\in \g^0$. 

Let $\wg=\g\oplus F$ be a non-trivial $1$-dimensional central 
extension of $\g$ defined by 
\begin{equation*}
[(x,\alpha), (y,\beta)] = ([x,y], \omega(x,y) ), \qquad x,y\in\g,
\end{equation*}
for a $2$-cocycle $\omega : \wedge^2 \g \to F$ defined as follows.
Let $r\in \h^0$ be a regular element, set 
$\omega(xt^i,\, yt^j)= \delta_{i+j,n+1}(x,[r,y])$, where $(,)$ is
the invariant scalar product in $\g$ (it is straightforward to
check that $\omega$  is a $2$-cocycle).

Suppose that $p$ is big enough so that we do not have to
divide by $p$ in the Campbell-Hausdorff formula
\begin{equation}
\exp(x) \exp(y) = \exp(x+y+\frac{1}{2}[x,y]+\dots), \qquad x,y\in \wg,
\end{equation}
where there are finitely many summands on the right hand side.
Then $\wG = \exp\wg = \{\exp(x) \mid x\in \wg\}$ is a finite
nilpotent group. Note that $\{\exp(\alpha) \mid \alpha\in F\}$
is a central subgroup of $\wG$ and that the latter is a central extension
of $G=\exp(\g)$.

The Cartan decomposition of $\g^0=\h^0\oplus \n^0_+\oplus \n^0_-$
yields the following decomposition of $\wg$ :
\begin{equation}
\wg = \wh\oplus \n_+\oplus \n_-,
\end{equation}
where $\wh = \h\oplus F$, $\h =\h^0 tF[t]/t^{n+1}\oplus F$, and 
$\n_\pm = \n^0_\pm tF[t]/t^{n+1}$.
For every $\lambda\in \h^*$ let us define a functional 
$\phi_\lambda\in \wg^*$
by setting
\begin{equation*}
\phi_\lambda|_{\n_+\oplus \n_-} = 0,\qquad
\phi_\lambda(h)  = \lambda(h),\qquad
\phi_\lambda(\alpha) = \alpha,
\end{equation*}
where $h\in\h,\, \alpha\in F$.
Let us denote $\wb_\pm = \n_\pm \oplus \wh$.

\begin{lemma}
\label{isotropic}
$\wb_\pm$ are Lie subalgebras of $\wg$ which are maximal isotropic subspaces
for the alternating bilinear form
\begin{equation}
b_\lambda(x, y) =  \phi_\lambda([x,y]), \qquad x,y\in \wg, \,
\lambda\in\h^*.
\end{equation}
\end{lemma}
\begin{proof}
It is clear from the above definitions that  $\wb_\pm$ are 
isotropic subspaces for  $b_\lambda$. Let us show that they are maximal.
The dimension $d$ of a maximal isotropic subspace of $b_\lambda$ is given by
\begin{equation}
d = \dim\wg - \frac{1}{2} \rank b_\lambda
  = \dim\wg - \frac{1}{2} \dim\O_\lambda,
\end{equation}
where $\O_\lambda$ is the orbit of the coadjoint action of $\wg$ containing
$\lambda$. We have $\dim\O_\lambda = \dim\wg - \dim\St(\phi_\lambda)$,
where
\begin{equation*}
\St(\phi_\lambda) = \{ x\in \wg \mid \Ad_x\phi_\lambda =0 \}
\end{equation*}
is the stabilizer of $\phi_\lambda$. Clearly, 
$\wh\subseteq  \St(\phi_\lambda)$.
In order to conclude that $d=\dim \wb_\pm$
it suffices to show that $\St(\phi_\lambda) =\wh$. 
Note that $\St(\phi_\lambda)$ is compatible with the Cartan decomposition
of $\wg$, so it is enough to prove that $\n_\pm \cap \St(\phi_\lambda) =0$.
Let $x\in \n_+\cap
\St(\phi_\lambda)$, where $x=x_mt^m+x_{m+1}t^{m+1}+\cdots$ with
$x_m\neq 0$. Let $y= bt^{n+1-m}$, where $b\in \n_-$ is such that
$(x_m,\, [r,b])\neq 0$. Then $\phi_\lambda([x,y]) = \omega(x_m, b)\neq 0$,
a contradiction that shows $\n_+ \cap \St(\phi_\lambda) =0$. One shows
that $\n_- \cap \St(\phi_\lambda) =0$ in a similar way.
\end{proof}

Let $\wB_\pm,\,N_\pm,\, H,\, \wH$ be the subgroups of $\wG$ corresponding
to Lie subalgebras $\wb_\pm,\,\n_\pm,\,\h,\,\wh$ of $\wg$, respectively.
Clearly, $\phi_\lambda$ gives rise
to the character $X_\lambda$ of $\wB_+$ (or $\wB_-$) :
\begin{equation}
X_\lambda(\exp(x)) = \psi\circ\phi_\lambda(x), \qquad x\in \wB_\pm,
\end{equation}
where $\psi : F \to k^\times$ is any fixed  non-trivial homomorphism.
The result of Kazhdan (\cite{K}, Proposition 2) together with 
Lemma~\ref{isotropic} imply that
\begin{equation}
V_\lambda := \Ind^{\wG}_{\wB_+} X_\lambda \cong \Ind^{\wG}_{\wB_-} X_\lambda 
\end{equation}
is an irreducible  representation of $\wG$ for every $\lambda\in \h^*$. 
Since $\exp{F}$ acts on $V_\lambda$ as a scalar defined by $\psi$, 
we have a family of projective  irreducible representations of $G$
parameterized by $\h^*$.
Note that since $\omega$ is trivial on $\h$ we have that 
$\wH \cong H \oplus \mathbb{Z}/p\mathbb{Z}$.

\begin{proposition}
$\{ V_\lambda \mid\lambda\in \h^*\} $ is a minimal dynamical datum for $G$.
%supported on $A=\exp(\sum_i \h t^i)$.
\end{proposition}
\begin{proof}
Let $\lambda,\mu\in \h^*$, $Y$ be a $G$-module,
and $k_\lambda$ denote a $1$-dimensional $\wH$-module determined by
$X_\lambda$, then we have
\begin{eqnarray*}
\Hom_{\wG}(V_\lambda\otimes V_\mu^*,\, Y)
&\cong& \Hom_{\wG}(\Ind^{\wG}_{\wB_+}X_\lambda \otimes
    \Ind^{\wG}_{\wB_-}X_{-\mu},\, Y) \\
&\cong& \Hom_{\wB_+}(k_\lambda\otimes \Ind^{\wG}_{\wB_-}X_{-\mu},\, Y) \\
&\cong& \Hom_{\wH}(k_\lambda\otimes k_{-\mu},\, Y) \cong Y[\lambda-\mu],
\end{eqnarray*}
where the second isomorphism is a consequence of the Frobenius reciprocity
and the third uses that $\Ind^{\wG}_{\wB_-}X_{\lambda}|_{\wB_+} \cong
\Ind_{\wH}^{\wB_+} X_{\lambda}$ (which holds since $\wG = \wB_+ N_-$,
$\wB_- = \wH N_-$, and $X_{\lambda}|_{N_-}=1$).  Therefore,
$V_\lambda\otimes V_\mu^*\cong \Ind_{\wH}^{\wG}\,X_{\lambda-\mu}
\cong  \Ind_{H}^{G}\,X_{\lambda-\mu}$, as required.
\end{proof}

\end{section}

%%%%%%%%%%%%%%%%%%%%%%%%%%%%%%%%%%%%%%%%%%%%%%%%%%%%%%%%%%%%%%%%%%%%%%%%%%%%%% 
%%%%%%%%%%%%%%%   BIBLIOGRAPHY %%%%%%%%%%%%%%%%%%%%%%%%%%%%%%%%%%%%%%%%%%%%%%% 
%%%%%%%%%%%%%%%%%%%%%%%%%%%%%%%%%%%%%%%%%%%%%%%%%%%%%%%%%%%%%%%%%%%%%%%%%%%%%  

\bibliographystyle{ams-alpha}
  
\end{document}